\documentclass{amsart}
\usepackage[T1]{fontenc}

\makeatletter

\providecommand{\LyX}{L\kern-.1667em\lower.25em\hbox{Y}\kern-.125emX\@}

\theoremstyle{plain}    
\newtheorem{thm}{Theorem}[section]
\numberwithin{equation}{section} 
\numberwithin{figure}{section} 
\theoremstyle{plain}    
\newtheorem{cor}[thm]{Corollary} 
\theoremstyle{plain}    
\newtheorem{lem}[thm]{Lemma} 
\theoremstyle{plain}    
\newtheorem{algorithm}[thm]{Algorithm} 
\theoremstyle{definition}
\newtheorem{defn}[thm]{Definition}
\theoremstyle{definition}
 \newtheorem{example}[thm]{Example}
\theoremstyle{remark}
\newtheorem{rem}[thm]{Remark}
\theoremstyle{remark}
\newtheorem*{rem*}{Remark}
\newenvironment{lyxlist}[1]
  {\begin{list}{}
    {\settowidth{\labelwidth}{#1}
     \setlength{\leftmargin}{\labelwidth}
     \addtolength{\leftmargin}{\labelsep}
     }}
  {\end{list}}

\makeatother

\begin{document}

\title{Definitive Computation of Bernstein-Sato Polynomials}

\author{Anton Leykin}

\begin{abstract}
Let \( n \) and \( d \) be positive integers, let \( k \) be a field and
let \( P(n,d;k) \) be the space of the polynomials in \( n \) variables of
degree at most \( d \) with coefficients in \( k \). Let \( B(n,d) \) be
the set of the Bernstein-Sato polynomials of all polynomials in \( P(n,d;k) \)
as \( k \) varies over all fields of characteristic \( 0 \). G. Lyubeznik
proved that \( B(n,d) \) is a finite set and asked if, for a fixed \( k \),
the set of the polynomials corresponding to each element of \( B(n,d) \) is
a constructible subset of \( P(n,d;k) \).

In this paper we give an affirmative answer to Lyubeznik's question by showing
that the set in question is indeed constructible and defined over \( \Bbb {Q} \),
i.e. its defining equations are the same for all fields \( k \). Moreover,
we construct an algorithm that for each pair \( (n,d) \) produces a complete
list of the elements of \( B(n,d) \) and, for each element of this list, an
explicit description of the constructible set of polynomials having this particular
Bernstein-Sato polynomial.
\end{abstract}
\maketitle

\section{Introduction}

Throughout this paper \( k \) is a field of characteristic \( 0 \), \( R_{n}(k)=k[x_{1},...,x_{n}] \)
is the ring of polynomials in \( n \) variables and \( A_{n}(k)=k\left\langle x_{1},...,x_{n},\partial _{1},...,\partial _{n}\right\rangle  \)
is the corresponding Weyl algebra, i.e. an associative \( k \)-algebra generated
by \( x \)'s and \( \partial  \)'s with the relations \( \partial _{i}x_{i}=x_{i}\partial _{i}+1 \)
for all \( i \). 

For every polynomial \( f \)\( \in R_{n}(k) \) there are \( b(s)\in k[s] \)
and \( Q(x,\partial ,s)\in A_{n}(k)[s] \) such that 
\begin{equation}
\label{func_eqn}
b(s)f^{s}=Q(x,\partial ,s)\cdot f^{s+1}.
\end{equation}
 The polynomials \( b(s) \) for which equation (\ref{func_eqn}) exists form
an ideal in \( k[s] \). The monic generator of this ideal is denoted by \( b_{f}(s) \)
and called the Bernstein-Sato polynomial of \( f. \) The roots of \( b_{f}(s) \)
are rational, in particular, \( b(s)\in \Bbb {Q}[s] \), where \( \Bbb {Q} \)
is the field of rational numbers. A good introduction to \( D \)-modules may
be found in \cite{bjork}.

The simplest characteristics of a polynomial \( f \) are its degree \( d \)
and its number of variables \( n \). This paper is motivated by the following
natural question: what can one say about \( b_{f}(s) \) in terms of \( n \)
and \( d \)? We give what may be regarded as a complete answer to this question.
Namely, we describe an algorithm that for fixed \( n \) and \( d \) gives
a complete list of all possible Bernstein-Sato polynomials and, for each polynomial
\( b(s) \) in this list, a complete description of the polynomials \( f \)
such that \( b_{f}(s)=b(s) \).

Let \( P(n,d;k) \) be the set of all the equivalence classes of the non-zero
polynomials of degree at most \( d \) in \( n \) variables with coefficients
in \( k \) modulo the equivalence relation \( f\sim g\, \Leftrightarrow \, f=c\cdot g\textrm{ } \)
for some \( 0\neq c\in k \). Note that \( b_{f}(s)=b_{g}(s) \) if \( f\sim g \).
We view \( P(n,d;k) \) as the set of the \( k \)-rational points of the projective
space \( \Bbb {P}(n,d;k)\cong \Bbb {P}_{k}^{N-1} \) where \( N \) is the number
of monomials in \( n \) variables of degree at most \( d \). G. Lyubeznik
\cite{lyubeznik} defined \( B(n,d) \) as the set of all the Bernstein-Sato
polynomials of all the polynomials from \( P(n,d;k) \) as \( k \) varies over
all fields of characteristic \( 0 \) and he proved that \( B(n,d) \) is a
finite set. He also asked if the subset of \( \Bbb {P}(n,d;k) \) corresponding
to a given element of \( B(n,d) \) is constructible. In this paper we give
an affirmative answer to Lyubeznik's question. The constructible sets in question
turn out to be definable over \( \Bbb {Q} \), i.e. their defining equations
and inequalities are the same for all fields \( k \).

A crucial ingredient in our proof is the fact, very recently discovered by T.
Oaku \cite{oaku} that there is an algorithm that, given a polynomial \( f \),
returns its Bernstein-Sato polynomial \( b_{f}(s) \). Using Oaku's algorithm
and our proof of the constructibility of the set of polynomials \( f \) having
a fixed \( b_{f}(s) \) we have developed an algorithm for the \emph{definitive
computation} of the Bernstein-Sato polynomials for each pair \( (n,d) \), by
which we mean that our algorithm, given \( n \) and \( d \), returns the list
of all the elements of \( B(n,d) \) and for each \( b(s)\in B(n,d) \), a finite
number of locally closed sets \( V_{i}=V_{i}'\setminus V_{i}'' \), where \( V_{i}' \)
and \( V_{i}'' \) are Zariski closed subsets of \( \Bbb {P}(n,d;\Bbb {Q}) \)
defined by explicit polynomial equations with rational coefficients, such that
for every field \( k \) of characteristic \( 0 \), the subset of \( P(n,d;k) \)
having \( b(s) \) as the Bernstein-Sato polynomial is the set of \( k \)-rational
points of 
\[
S(b(s),k)=(\cup _{i}V_{i})\otimes _{\Bbb {Q}}k\subset \Bbb {P}(n,d;\Bbb {Q})\otimes _{\Bbb {Q}}k=\Bbb {P}(n,d;k).\]

The definitive computation for fixed \( n \) and \( d \) could be useful in
a number of ways. For example, it would produce an algorithm for the computation
of \( b_{f}(s) \) for \( f\in P(n,d;k) \) that is likely to be considerably
more efficient than all other currently available algorithms. It would also
produce the smallest integer \( t \) such that \( R_{n}(k)_{f} \) is generated
by \( \frac{1}{f^{t}} \) as an \( A_{n}(k) \)-module for all \( f \) of degree
at most \( d \) (this integer is denoted \( t(n,d) \) in \cite{lyubeznik}).
Moreover, using a similar technique we develop an algorithm for a \emph{quasi-definitive
computation} of the annihilator of \( \frac{1}{f^{t}} \) in \( A_{n}(k) \),
which, provided \( t \) is known, gives a presentation of \( R_{n}(k)_{f} \)
as an \( A_{n}(k) \)-module (see Example \ref{exaAnnFs}). We call it quasi-definitive
because its output is not uniquely determined (see Remark \ref{rem_quasidefinitive}).
The last algorithm is particularly important for U. Walther's algorithmic computation
of local cohomology modules \cite{Uli}. These applications are discussed in
the next section. More applications will undoubtedly arise in the future. 

The results of this paper are a part of my thesis. I would like to thank my
advisor Gennady Lyubeznik for suggesting this problem to me.

\section{Examples and Discussion}

Our algorithms have been implemented as scripts written in the Macaulay 2 programming
language (see \cite{M2}). In this section we give some examples of actual computations
and discuss possible uses of the results of computation.

\begin{example}
\label{exB22}If \( n=2 \) and \( d=2 \) then 
\[
f=a_{20}x^{2}+a_{11}xy+a_{02}y^{2}+a_{10}x+a_{01}y+a_{00},\]
so \( P(2,2;k) \) is the set of the k-rational points of the projective space
\( \Bbb {P}(2,2;k)=\Bbb {P}^{5}_{k} \) with the homogeneous coordinate ring
\( k[a_{ij}] \), \( i,j=0,1,2 \). It takes our program less than 20 minutes
on 300MHz Pentium-II machine to produce 
\[
B(2,2)=\{1,\, s+1,\, (s+1)^{2},\, (s+1)(s+\frac{1}{2})\}\]
 and give a description of the corresponding constructible sets of polynomials
from \( B(2,2) \) which is essentially equivalent to the following:

\( \bullet  \) \( b_{f}(s)=1 \) iff \( f\in V_{1}=V_{1}'\setminus V_{1}'' \),
where \( V_{1}'=V(a_{1,1},a_{0,1},a_{0,2},a_{1,0},a_{2,0}) \), while \( V_{1}''=V(a_{0,0}) \),

\( \bullet  \) \( b_{f}(s)=s+1 \) iff \( f\in V_{2}=(V_{2}'\setminus V_{2}'')\cup (V_{3}'\setminus V_{3}'') \),
where \( V_{2}'=V(0) \), \( V_{2}''=V(\gamma _{1}) \), \( V_{3}'=V\left( \gamma _{2},\gamma _{3},\gamma _{4}\right)  \),
\( V_{3}''=V\left( \gamma _{3},\gamma _{4},\gamma _{5},\gamma _{6},\gamma _{7},\gamma _{8}\right)  \), 

\( \bullet  \) \( b_{f}(s)=(s+1)^{2} \) iff \( f\in V_{4}'\setminus V_{4}'' \),
where \( V_{4}'=V(\gamma _{1}) \), \( V_{4}''=V\left( \gamma _{2},\gamma _{3},\gamma _{4}\right)  \),

\( \bullet  \) \( b_{f}(s)=(s+1)(s+\frac{1}{2}) \) iff \( f\in V_{5}'\setminus V_{5}'' \),
where \( V_{5}'=V\left( \gamma _{3},\gamma _{4},\gamma _{5},\gamma _{6},\gamma _{7},\gamma _{8}\right)  \),
while \( V_{5}''=V(a_{1,1},a_{0,1},a_{0,2},a_{1,0},a_{2,0}) \),

where \( \gamma _{i} \) may be looked up in this list:

\( \gamma _{1}=a_{0,2}a_{1,0}^{2}-a_{0,1}a_{1,0}a_{1,1}+a_{0,0}a_{1,1}^{2}+a_{0,1}^{2}a_{2,0}-4a_{0,0}a_{0,2}a_{2,0} \),

\( \gamma _{2}=2a_{0,2}a_{1,0}-a_{0,1}a_{1,1} \),

\( \gamma _{3}=a_{1,0}a_{1,1}-2a_{0,1}a_{2,0} \), 

\( \gamma _{4}=a_{1,1}^{2}-4a_{0,2}a_{2,0} \),

\( \gamma _{5}=2a_{0,2}a_{1,0}-a_{0,1}a_{1,1} \),

\( \gamma _{6}=a_{0,1}^{2}-4a_{0,0}a_{0,2} \),

\( \gamma _{7}=a_{0,1}a_{1,0}-2a_{0,0}a_{1,1} \), 

\( \gamma _{8}=a_{1,0}^{2}-4a_{0,0}a_{2,0} \).

It is not hard to see that this definitive computation agrees with the well-known
result that \( b_{f}(s)=1 \) iff \( f \) is constant, \( b_{f}(s)=s+1 \)
iff \( f \) is non-constant and non-singular, and \( b_{f}(s)=(s+1)^{2} \)
( resp. \( b_{f}(s)=(s+1)(s+\frac{1}{2}) \) ) iff f can be reduced to \( xy \)
(resp. \( x^{2} \)) by a linear change of variables.

\hfill{}
\end{example}
The definitive computation for fixed \( n \) and \( d \) is likely to lead
to a considerably more efficient way of computing \( b_{f}(s) \) for \( f\in P(n,d;k) \).
Namely, to compute \( b_{f}(s) \) for a concrete polynomial \( f \) one just
has to ``search the database'', i.e. check which of the constructible sets
this polynomial belongs to. Since there are finitely many of them and each one
is described by explicit equations and inequalities in the coefficients of \( f \)
and each \( f \) belongs to a unique one, we get a straightforward algorithm
for computing \( b_{f}(s) \) for all \( f\in P(n,d;k) \).

All other known algorithms for computing \( b_{f}(s) \) involve Gröbner bases
computations. Often \( b_{f}(s) \) is not very big but its computation is enormous
because of the \char`\"{}intermediate explosion\char`\"{} caused by the fact
that Gröbner bases computations are very time and memory consuming. But the
algorithm of ``searching the database'' does not involve any Gröbner bases
at all! For this reason it is likely to be considerably more efficient in computing
\( f \) for \( f\in P(n,d;k) \), especially if the field \( k \) is the fraction
field of some finitely generated \( \Bbb {Q} \)-algebra, so that ordinary arithmetic
operations in \( k \) and hence Gröbner bases computations are especially expensive.

Certainly the algorithm just described requires ``setting up the database''.
A definitive computation for \( n \) and \( d \) must be performed just once.
This part may be done on a ``powerful computer'' (we have in mind implementing
some parallel processing techniques) and the results of this computation may
then be stored in a file accessible for ``not-so-powerful'' machines, which
are capable of performing the ``search the database'' part. However a definitive
computation even for rather small values of \( n \) and \( d \) with the modest
computer resources at our disposal and with the current level of efficiency
of our program faces its own \char`\"{}intermediate explosion\char`\"{} problem.

\hfill{}

\begin{example}
If \( n=2 \) and \( d=3 \) then 
\begin{eqnarray*}
f & = & a_{3,0}x^{3}+a_{2,1}x^{2}y+a_{1,2}xy^{2}+a_{0,3}y^{3}\\
 & + & a_{2,0}x^{2}+a_{1,1}xy+a_{0,2}y^{2}+a_{1,0}x+a_{0,1}y+a_{0,0},
\end{eqnarray*}
so \( P(2,3;k) \) is the set of the k-rational points of \( \Bbb {P}(2,3;k)=\Bbb {P}_{k}^{9} \)
with the homogeneous coordinate ring that involves \( 10 \) variables. Our
program exhausts all available memory, 128Mb, of the computer after about 3
hours and stops without producing an answer. However, a somewhat creative use
of our program enables us to give a complete list of all the elements of \( B(2,3) \)
(but not the explicit descriptions of the constructible sets corresponding to
each element of \( B(2,3) \)): \emph{}
\end{example}
Since for any nonsingular polynomial its Bernstein-Sato polynomial is equal
to \( s+1 \), it remains to consider the case where our \( f\in \Bbb {P}(2,3;k) \)
possesses a singularity at some point \( (x_{0},y_{0}) \). Keeping in mind
that the \emph{}Bernstein-Sato polynomial is stable under any linear substitution
of variables, we may get rid of its linear part via the substitution \( x\mapsto x-x_{0} \),
\( y\mapsto y-y_{0} \), i.e. \( f \) takes the form
\[
f=(ax^{3}+bx^{2}y+cxy^{2}+dy^{3})+(a'x^{2}+b'xy+c'y^{2}).\]
 Now it is easy to see that by homogeneous linear transformation the quadratic
part may be shaped to one of the forms \( 0 \), \( xy \), \( x^{2} \). Therefore
it is enough to compute the Bernstein-Sato polynomial for the following polynomials:
\begin{eqnarray*}
f_{1} & = & ax^{3}+bx^{2}y+cxy^{2}+dy^{3},\\
f_{2} & = & (ax^{3}+bx^{2}y+cxy^{2}+dy^{3})+xy,\\
f_{3} & = & (ax^{3}+bx^{2}y+cxy^{2}+dy^{3})+x^{2}.
\end{eqnarray*}
 Our program returns the complete sets of possible Bernstein-Sato polynomials
for \( f_{1} \) in 22 minutes, for \( f_{2} \) in 16 minutes and for \( f_{3} \)
in 21 minutes. Of course, in each of the three cases our program produces an
explicit description of the corresponding constructible set in \( \Bbb {P}^{3}_{k} \)
(each of \( f_{i} \) contains 4 indeterminate coefficients) for each element
\( b(s)\in B_{f_{i}} \). We omit these and list only the Bernstein-Sato polynomials:

\begin{eqnarray*}
B_{f_{1}} & = & \{\, (s+1)^{2}(s+\frac{2}{3})(s+\frac{4}{3}),\\
 &  & (s+1)^{2}(s+\frac{1}{2}),\\
 &  & (s+1)(s+\frac{2}{3})(s+\frac{1}{3})\, \};\\
B_{f_{2}} & = & \{\, (s+1)^{2}\, \};\\
B_{f_{3}} & = & \{\, (s+1)(s+\frac{7}{6})(s+\frac{5}{6}),\\
 &  & (s+1)^{2}(s+\frac{3}{4})(s+\frac{5}{4}),\\
 &  & (s+1)^{2}(s+\frac{1}{2}),\\
 &  & (s+1)(s+\frac{1}{2})\, \}.
\end{eqnarray*}

Thus

\begin{eqnarray*}
B(2,3) & = & \{\, (s+1)^{2}(s+\frac{2}{3})(s+\frac{4}{3}),\\
 &  & (s+1)^{2}(s+\frac{1}{2}),\\
 &  & (s+1)(s+\frac{2}{3})(s+\frac{1}{3}),\\
 &  & (s+1)^{2},\\
 &  & (s+1)(s+\frac{7}{6})(s+\frac{5}{6}),\\
 &  & (s+1)^{2}(s+\frac{3}{4})(s+\frac{5}{4}),\\
 &  & (s+1)(s+\frac{1}{2}),\\
 &  & s+1,\\
 &  & 1\, \}.
\end{eqnarray*}

As was mentioned above, only the efficiency of the algorithm and the current
efficiency of computer hardware and software obstruct us from getting a complete
description of the constructible sets that correspond to the polynomials above. 

As was pointed out in \cite{lyubeznik}, \( t(n,d) \) (which is defined in
the last paragraph of the preceding section) is the largest absolute value of
all the negative integer roots of all the polynomials in \( B(n,d) \). Thus
we get 

\begin{cor}
\( t(2,3)=1 \), i.e. if \( f\in R_{2}(k) \) is of degree at most \( 3 \),
then \( \frac{1}{f} \) generates \( R_{2}(k)_{f} \) as an \( A_{2}(k) \)-module.
\end{cor}
To compute the localization of \( R_{n}(k) \) at a polynomial \( f\neq 0 \)
one needs to compute \( \text {Ann}f^{s}\subset A_{n}(k)[s] \) and take \( N=\left. \text {Ann}f^{s}\right| _{s=a} \),
where \( a \) is the minimal integer root of \( b_{f}(s) \). Then \( R_{n}(k)_{f}=A_{n}(k)/N \)
as an \( A_{n}(k) \)-module (see Section \ref{secOaku} below). 

Using a technique similar to that for computing Bernstein-Sato polynomials,
we constructed an algorithm for a quasi-definitive computation of \( \text {Ann}f^{s} \)
for all \( f\in P(n,d;k) \). By this we mean an explicit subdivision of \( \Bbb {P}(n,d;k) \)
into a finite union of constructible subsets and for each such subset \( V, \)
an explicit finite set of elements \( \beta _{1},\beta _{2},...\in A_{n}(k)[a_{i_{1}...i_{n}}][s] \)
with \( i_{1}+...+i_{n}\leq d \), such that \( \text {Ann}(f^{s})=(\beta _{1}',\beta _{2}',...) \)
for every \( f\in V \), where \( \beta _{i}' \) is the image of \( \beta _{i} \)
under the specialization of the \( a_{i_{1}...i_{n}} \) to the corresponding
coefficients of \( f \). 

\begin{example}
\label{exaAnnFs}Here is what we got for \( P(2,2;k) \) (See Example \ref{exB22}
for notation):

\( \bullet  \) \( \text {Ann}(f^{s})=(\beta _{1},\beta _{2},\beta _{3}) \)
if \( f\in (V_{1}'\setminus V_{1}'')\cup (V_{2}'\setminus (V_{2,1}''\cup V_{2,2}'')) \),
where \( V_{1}'=V(0) \), \( V_{1}''=V(\gamma _{1}) \), \( V_{2}'=V(\gamma _{2},\gamma _{3},\gamma _{4}) \),
\( V_{2,1}''=V(a_{1,1},\, a_{0,2},\, a_{0,1}) \) and \( V_{2,2}''=V(\gamma _{2},\gamma _{3},\gamma _{4},\gamma _{5},\gamma _{6},\gamma _{7}) \);

\( \bullet  \) \( \text {Ann}(f^{s})=(\beta _{1},\beta _{4}) \) if \( f\in V_{3}'\setminus V_{3}'' \),
where \( V_{3}'=V(\gamma _{1}) \), while \( V_{3}''=V(\gamma _{2},\gamma _{3},\gamma _{4}) \);

\( \bullet  \) \( \text {Ann}(f^{s})=(\beta _{5},\beta _{6}) \) if \( f\in V_{4}'\setminus (V_{4,1}''\cup V_{4,2}'') \),
where \( V_{4}'=V(\gamma _{2},\gamma _{3},\gamma _{4},\gamma _{5},\gamma _{6},\gamma _{7}) \),
\( V_{4,1}''=V(a_{1,0},\, a_{2,0},\, a_{1,1},\, \gamma _{5}) \) and \( V_{4,2}''=(a_{1,1},\, a_{0,2},\, a_{0,1},\, \gamma _{7}) \);

\( \bullet  \) \( \text {Ann}(f^{s})=(\beta _{7},\beta _{8}) \) if \( f\in V_{5}'\setminus V_{5}'' \),
where \( V_{5}'=V(a_{1,0},\, a_{2,0},\, a_{1,1},\, \gamma _{5}) \), while \( V_{5}''=V(a_{1,1},\, a_{0,1},\, a_{0,2},\, a_{1,0},\, a_{2,0}) \); 

\( \bullet  \) \( \text {Ann}(f^{s})=(\beta _{9},\beta _{10}) \) if \( f\in V_{6}'\setminus V_{6}'' \),
where \( V_{6}'=V(a_{1,1},\, a_{0,2},\, a_{0,1},\, \gamma _{7}) \), while \( V_{5}''=V(a_{1,1},\, a_{0,1},\, a_{0,2},\, a_{1,0},\, a_{2,0}) \); 

\( \bullet  \) \( \text {Ann}(f^{s})=(\beta _{9},\beta _{11}) \) if \( f\in (V_{7}'\setminus V_{7}'')\cup V_{8}' \),
where \( V_{7}'=V(a_{1,1},\, a_{0,2},\, a_{0,1}) \), \( V_{7}''=V(a_{1,1},\, a_{0,2},\, a_{0,1},\, \gamma _{7}) \)
and \( V_{8}'=V(a_{1,1},\, a_{0,1},\, a_{0,2},\, a_{1,0},\, a_{2,0}) \);

where the polynomials \( \beta _{i} \) are listed below:

\( \beta _{1}=a_{1,1}x_{1}dx_{1}+2a_{0,2}x_{2}dx_{1}-2a_{2,0}x_{1}dx_{2}-a_{1,1}x_{2}dx_{2}+a_{0,1}dx_{1}-a_{1,0}dx_{2} \),

\( \beta _{2}=a_{1,1}a_{2,0}x_{1}^{2}dx_{1}+a_{1,1}^{2}x_{1}x_{2}dx_{1}+a_{0,2}a_{1,1}x_{2}^{2}dx_{1}-2a_{2,0}^{2}x_{1}^{2}dx_{2}-2a_{1,1}a_{2,0}x_{1}x_{2}dx_{2} \)

\( \, \, -2a_{0,2}a_{2,0}x_{2}^{2}dx_{2}-a_{1,1}^{2}sx_{2}+4a_{0,2}a_{2,0}sx_{2}+a_{1,0}a_{1,1}x_{1}dx_{1}+a_{0,1}a_{1,1}x_{2}dx_{1} \)

\( \, \, -2a_{1,0}a_{2,0}x_{1}dx_{2}-2a_{0,1}a_{2,0}x_{2}dx_{2}-a_{1,0}a_{1,1}s+2a_{0,1}a_{2,0}s+a_{0,0}a_{1,1}dx_{1}-2a_{0,0}a_{2,0}dx_{2} \),

\( \beta _{3}=a_{2,0}x_{1}^{2}dx_{2}+a_{1,1}x_{1}x_{2}dx_{2}+a_{0,2}x_{2}^{2}dx_{2}-a_{1,1}sx_{1}-2a_{0,2}sx_{2}+a_{1,0}x_{1}dx_{2} \)

\( \, \, +a_{0,1}x_{2}dx_{2}-a_{0,1}s+a_{0,0}dx_{2} \),

\( \beta _{4}=a_{1,1}^{2}x_{1}dx_{1}-4a_{0,2}a_{2,0}x_{1}dx_{1}+a_{1,1}^{2}x_{2}dx_{2}-4a_{0,2}a_{2,0}x_{2}dx_{2}-2a_{1,1}^{2}s \)

\( \, \, +8a_{0,2}a_{2,0}s-2a_{0,2}a_{1,0}dx_{1}+a_{0,1}a_{1,1}dx_{1}+a_{1,0}a_{1,1}dx_{2}-2a_{0,1}a_{2,0}dx_{2} \),

\( \beta _{5}=a_{1,1}dx_{1}-2a_{2,0}dx_{2} \),

\( \beta _{6}=2a_{2,0}x_{1}dx_{2}+a_{1,1}x_{2}dx_{2}-2a_{1,1}s+a_{1,0}dx_{2} \),

\( \beta _{7}=dx_{1} \),

\( \beta _{8}=2a_{0,2}x_{2}dx_{2}-4a_{0,2}s+a_{0,1}dx_{2} \),

\( \beta _{9}=dx_{2} \),

\( \beta _{10}=2a_{2,0}x_{1}dx_{1}-4a_{2,0}s+a_{1,0}dx_{1} \),

\( \beta _{11}=a_{2,0}x_{1}^{2}dx_{1}-2a_{2,0}sx_{1}+a_{1,0}x_{1}dx_{1}-a_{1,0}s+a_{0,0}dx_{1} \),

and the polynomials \( \gamma _{i} \) are in this list:

\( \gamma _{1}=a_{0,2}a_{1,0}^{2}-a_{0,1}a_{1,0}a_{1,1}+a_{0,0}a_{1,1}^{2}+a_{0,1}^{2}a_{2,0}-4a_{0,0}a_{0,2}a_{2,0} \),

\( \gamma _{2}=2a_{0,2}a_{1,0}-a_{0,1}a_{1,1} \),

\( \gamma _{3}=a_{1,0}a_{1,1}-2a_{0,1}a_{2,0} \), 

\( \gamma _{4}=a_{1,1}^{2}-4a_{0,2}a_{2,0} \),

\( \gamma _{5}=a_{0,1}^{2}-4a_{0,0}a_{0,2} \), 

\( \gamma _{6}=a_{0,1}a_{1,0}-2a_{0,0}a_{1,1} \), 

\( \gamma _{7}=a_{1,0}^{2}-4a_{0,0}a_{2,0} \).
\end{example}

\section{Constructible Sets}

In this section we describe some of the properties of constructible sets that
are used in the next section. We recall that a set is constructible iff it is
a finite union of locally closed sets and a set is locally closed iff it is
the difference of two closed sets.

\begin{thm}
Let \( C \) be a constructible subset of a variety \( X \). Then \( C \)
may be presented uniquely as a disjoint union \( \bigcup ^{m}_{i=1}(V_{i}'\setminus V_{i}'') \),
where for all \( i \) the sets \( V_{i}' \) and \( V_{i}'' \) are closed,
\( V_{1}'\supset V_{1}''\supset V_{2}'\supset V_{2}''\supset ...\supset V_{m}'\supset V_{m}'' \)
and no two consequent sets in this chain have an irreducible component in common.
We call it a canonical presentation of \( C \) as a union of locally closed
subsets.
\end{thm}
\begin{proof}
Let \( d(C) \) be the maximal dimension of an irreducible component in \( \bar{C} \).
The only possible choices for \( V_{1}' \) and \( V_{2}' \). Now \( V_{1}'=\bar{C} \)
and \( V_{1}''=\overline{V_{1}'\setminus C} \) and let \( C_{1}=C\cap V_{1}'' \).
Note that \( d(C_{1})<d(C) \) and we may assume by induction on \( d \) that
the chain \( V_{2}'\supset V_{2}''\supset ...\supset V_{m}'\supset V_{m}'' \)
such that \( C'=\bigcup ^{m}_{i=2}(V_{i}'\setminus V_{i}'') \) exists and is
unique. Then \( V_{1}'\supset V_{1}''\supset V_{2}'\supset V_{2}''\supset ...\supset V_{m}'\supset V_{m}'' \)
is the unique chain for \( C \), which satisfies the condition in the statement.
\end{proof}
\begin{rem}
\label{remark_constr_presentation}There is an algorithmic way for constructing
such a presentation, starting with \( C \) presented as a union of nonempty
sets \( W_{\alpha }\setminus (W_{\alpha }^{(1)}\cup ...\cup W_{\alpha }^{(h_{\alpha })}) \),
where \( W_{\alpha } \) and \( W_{\alpha }^{(i)} \) are closed irreducible
subsets and \( W_{\alpha }\supset W_{\alpha }^{(i)} \) for all \( i \). Let
\( d(C)=\max _{\alpha }\dim W_{\alpha } \) (which agrees with the definition
in the proof of the theorem). 

Let \( V_{1}' \) be the union of all maximal elements in the set \( \{W_{\alpha }\} \)
and \( V_{1}'' \) be the union of all \( W_{\alpha }^{(i)} \) that are minimal
with the following property: there is a set of pairs \( \{(\alpha _{j},i_{j})\}_{j=1}^{l} \)
such that \( W_{\alpha _{1}} \) is a component of \( V_{1}' \), \( W_{\alpha _{l}}^{(i_{l})}=W_{\alpha }^{(i)} \)
and \( W_{\alpha _{j}}^{(i_{j})}\supset W_{\alpha _{j}-1} \) for all \( j=2,...,l \)
. Now \( d(C\setminus (V_{1}'\setminus V_{1}'')) \) is less than \( d(C) \),
therefore, we may assume again by induction on \( d \) that we are able to
construct the rest of \( V_{i}' \) and \( V_{i}'' \). 
\end{rem}
\begin{lem}
\label{lemma_constructible}Let \( X \) be a variety and \( f:X\rightarrow Y \)
a map into any finite set \( Y \). Then \( f^{-1}(y) \) is constructible for
every \( y\in Y \) iff for every closed irreducible subvariety \( X'\subset X \)
there is an open \( U\subset X' \) such that \( \left. f\right| _{U} \) is
a constant function.
\end{lem}
\begin{proof}
Assume the second part holds. Take any \( y\in Y \) and let \( Z=f^{-1}(y) \).
Let \( n=\dim X \) and assume the lemma is proved for dimensions less then
\( n \). First of all, since \( X \) is a finite union of its irreducible
components, we may proceed assuming that \( X \) is irreducible. Let \( U \)
be an open subset of \( X \) such that \( f(u)=y' \) for all \( u\in U \).
If \( y'\neq y \) then \( Z\subset X\setminus U \), which has dimension less
than \( n \) and, therefore, \( Z \) is constructible by the induction assumption.
If \( y=y' \) then \( (Z\setminus U)\subset (X\setminus U) \) is constructible,
hence so is \( Z=U\cup (Z\setminus U) \). 

It remains to check the case \( \dim X=0 \), in which \( X \) is a finite
set of points and is certainly constructible. 

Conversely, assume that \( f^{-1}(y) \) is constructible for every \( y\in Y \).
Let \( X'\subset X \) be a closed irreducible subvariety. Then \( X'=\bigcup _{y\in Y}(f^{-1}(y)\cap X') \)
and, since \( Y \) is a finite set and \( X' \) is irreducible, the closure
of \( X_{y}'=f^{-1}(y)\cap X' \) for some \( y\in Y \) is equal to \( X' \).
But \( X_{y}' \) is constructible, hence it is open in its closure \( \overline{X_{y}'}=X' \).
\end{proof}

\section{Parametric Gröbner Bases}

This section describes an approach to computing parametric Gröbner bases in
Weyl algebras. A good source on computing Gröbner bases in non-commutative algebras
is \cite{non_comm_GBs}. For a discussion of parametric Gröbner bases, which
leads to the notion of comprehensive Gröbner bases, see \cite{comprehensive}
for the commutative case and \cite{paramGB} for the case of solvable algebras.
However, everything that is needed for this paper is stated and proved in this
section.

Let \( C=k[\bar{a}] \) (\( \bar{a}=\{a_{1},...,a_{m}\} \)) be the ring of
parameters and \( R=C\left\langle \bar{y},\bar{x},\bar{\partial }\right\rangle  \)
be the ring of non-commutative polynomials in \( \bar{y}=\{y_{1},...,y_{l}\} \),
\( \bar{x}=\{x_{1},...,x_{n}\} \) and \( \bar{\partial }=\{\partial _{1},...,\partial _{n}\} \)
with coefficients in \( C \), where \( \bar{x} \) and \( \bar{\partial } \)
satisfy the same relations as in a Weyl algebra and \( \bar{y} \) is contained
in the center of \( R \).

\begin{defn}
\label{defSpecialization} For a prime \( P \) in \( C \), we shall call the
natural map \( C\rightarrow k(P) \) as well as the induced map \( R=C\left\langle \bar{y},\bar{x},\bar{\partial }\right\rangle \rightarrow k(P)\left\langle \bar{y},\bar{x},\bar{\partial }\right\rangle  \),
where \( k(P) \) is the residue field at \( P \), the \emph{specialization}
at the point \( P \) and denote both maps by \( \sigma _{P} \).
\end{defn}
The next result is similar to Oaku's Proposition 7 in \cite{oaku}.

Let \( < \) be an order on monomials in \( a \), \( y \), \( x \) and \( \partial  \)
such that every \( a_{i} \) is \( << \) than any of \( x_{j} \), \( y_{j} \)
or \( \partial _{j} \) (i.e. the order \( < \) eliminates \( x_{j} \), \( y_{j} \)
and \( \partial _{j} \)). Assume \( G \) is a finite Gröbner basis in \( R \),
then we claim that \( \sigma _{P}(G)=\{\sigma _{P}(g)\, |\, g\in G\} \) is
a Gröbner basis in \( \sigma _{P}(R) \) for ``almost'' every \( P\in \text {Spec}\, C \).
Namely,

\begin{lem}
For any \( G \) there exists a polynomial \( h\in C \) such that \( \sigma _{P}(G) \)
is a Gröbner basis for every \( P \) not containing \( h \). 
\end{lem}
\begin{proof}
We need to make some definitions. For a polynomial \( f \) let \( inM(f) \)
be the initial monomial \( inC(f) \) the initial coefficient such that \( in(f)=inC(f)\cdot inM(f) \)
the initial term of \( f \). Also for \( f\in R \) let \( inM_{*}(f)\in \left\langle \bar{y},\bar{x},\bar{\partial }\right\rangle  \)
and \( inC_{*}(f)\in C \) be the initial monomial and the initial coefficient
of \( f \) viewed as a polynomial in \( x,y,\partial  \) with coefficients
in \( C \) with respect to \( \prec  \), the restriction of \( < \) to \( \left\langle \bar{y},\bar{x},\bar{\partial }\right\rangle  \). 

One obvious observation is that a specialization \( \sigma _{P}:(R,<)\rightarrow (\sigma _{P}(R),\prec ) \)
preserves the order.

Let \( h=\prod _{g\in G}inC_{*}(g)\in C \). Consider any \( P\in \text {Spec}\, C \)
not containing \( h \). Take a polynomial \( f' \) in ideal of \( \sigma _{P}(R) \)
generated by \( \sigma _{P}(G) \), then there is \( f \) such that \( f'=\sigma _{P}(f) \)
and \( inM_{*}(f)=inM(f') \). Since \( G \) is a Gröbner basis in \( R \),
we have \( inM(g)|inM(f) \) for some \( g\in G \), which means that \( inM_{*}(g)|inM_{*}(f) \).
Now, \( inM(\sigma _{P}(g))=inM_{*}(g) \), because \( inC_{*}(g)\notin P \).
Thus \( inM(\sigma _{P}(g))|inM(\sigma _{P}(f)) \), which proves that \( \sigma _{P}(G) \)
is a Gröbner basis.
\end{proof}
\begin{rem}
The statement of the lemma is true for reduced Gröbner bases as well.
\end{rem}
The lemma leads to the following 

\begin{algorithm}
\label{alg_param} 
\end{algorithm}
\vspace{0.3cm}
{\centering \begin{tabular}{lll}
Input: &
\( F' \): &
a finite set of generators for a prime ideal \( Q\subset C \). \\
&
\( F \): &
a finite set of generators of a left ideal \( I\subset R \) containing \( QA_{n} \),
\\
Output: &
\( G \): &
a (reduced) Gröbner basis in \( R \) with respect to \( < \), \\
&
\( h \): &
a polynomial in \( C \), which we shall call an \emph{exceptional polynomial},
\\
&
&
such that for any \( P\in \text {Spec}(k[a_{1},...,a_{m}]) \), \( P\supset Q \)
and \( h\notin P \)\\
&
&
the ideal \( \sigma _{P}(I)\subset \sigma _{P}(R) \) has a \( \sigma _{P}(G) \)
as a (reduced) Gröbner \\
&
&
basis with respect to \( \prec  \).\\
\end{tabular}\par}
\vspace{0.3cm}

\begin{enumerate}
\item Compute a Gröbner basis \( G \) of \( I+QR \) (which is generated by \( F\cup F' \))
.
\item Return \( G \) and \( h=\prod _{g\in G\setminus Q}inC_{*}(g) \).
\end{enumerate}
\begin{rem}
If all polynomials in \( F' \) and all \( C \)-coefficients of all elements
of \( F' \) are homogeneous, then so is the exceptional polynomial \( h \).
\end{rem}

\section{Oaku's Algorithm\label{secOaku}}

The original algorithm of T.Oaku for computing the Bernstein-Sato polynomial
appeared in \cite{oaku}. However there exist several modifications of the algorithm
(see \cite{SST} for example). For our needs a version of the algorithm described
in \cite{Uli} will be utilized. 

Let \( f\in R_{n}(k) \). Denote by \( \text {Ann}f^{s} \) the ideal of all
elements in \( A_{n}(k)[s] \) annihilating \( f^{s} \). The following algorithm
is Algorithm 4.4. from \cite{Uli} with \( L=(\partial _{1},...,\partial _{n}) \).

\begin{algorithm}
\label{AnnFs}
\end{algorithm}
\vspace{0.3cm}
{\centering \begin{tabular}{lll}
Input: &
\( f \):&
a polynomial in \( R_{n}(k) \) ,\\
Output: &
\( \{P_{j}'\} \):&
generators of \( \text {Ann}f^{s} \)\\
\end{tabular}\par}
\vspace{0.3cm}

\begin{enumerate}
\item Set \( Q=\{\partial _{i}+\frac{df}{dx_{i}}\partial _{t},t\} \).
\item Homogenize all \( q_{i}\in Q \) using the new variable \( y_{1} \) with respect
to the weight \( w \), where \( w(t)=w(y_{1})=1 \), \( w(\partial _{t})=w(y_{2})=-1 \),
\( w(x_{i})=w(\partial _{i})=0 \). Denote the homogenized elements \( q_{i}^{h} \).
\item Compute a Gröbner basis for the ideal generated by \( q_{1}^{h},\, ...\, ,\, q_{r}^{h},\, 1-y_{1}y_{2} \)
in \( A_{n+1}[y_{1},y_{2}] \) with respect to an order eliminating \( y_{1},y_{2} \).
\item Select the operators \( \{p_{j}\}_{1}^{b} \) in this basis which do not contain
\( y_{1},y_{2} \).
\item For each \( p_{j} \), if \( w(p_{j})>0 \) then replace \( p_{j} \) by \( p_{j}'=\partial _{t}^{w(P_{j})}p_{j} \)
else replace \( p_{j} \) by \( p_{j}'=t^{-w(P_{j})}p_{j} \).
\item Return the operators \( \{p_{j}'\}_{1}^{b} \).
\end{enumerate}
The following is Algorithm 4.6 in \cite{Uli}.

\begin{algorithm}
\label{oakuAlg}
\end{algorithm}
\vspace{0.3cm}
{\centering \begin{tabular}{lll}
Input:&
\( f: \)&
o polynomial in \( R_{n}(k) \),\\
Output:&
\( b_{f}(s) \)&
the Bernstein-Sato polynomial of \( f \).\\
\end{tabular}\par}
\vspace{0.3cm}

\begin{enumerate}
\item Determine \( \text {Ann}f^{s} \) following Algorithm \ref{AnnFs}.
\item Find a reduced Gröbner basis for the ideal \( \text {Ann}f^{s}+A_{n}[s]\cdot f \)
using an order that eliminates \( x \) and \( \partial  \).
\item Return the unique element in the basis contained in \( k[s] \).
\end{enumerate}

\section{The Main Results\label{sec_algorithm}}

Consider \( \Bbb {P}(n,d;k) \) with the coordinate ring \( C=k[\bar{a}] \),
where \( \bar{a}=\{a_{\alpha }\, :\, |\alpha |\leq d\} \). Let \( f=\sum _{|\alpha |\leq d}a_{\alpha }x^{\alpha } \). 

\begin{defn}
\label{SBSdef}Let \( b(s)\in B(n,d) \). The corresponding set \( S(b(s),k)\subset \Bbb {P}(n,d;k) \)
in \( \Bbb {A}_{k}^{N} \) is defined as the set of all the points \( P\in \Bbb {P}(n,d;k) \)
such that \( b_{\sigma _{P}(f)}(s)=b(s) \). (We view points in \( \Bbb {P}(n,d;k) \)
as homogeneous primes in \( C \). See Definition \ref{defSpecialization} for
\( \sigma _{P}(f) \).)
\end{defn}
\hfill{}

Let \( Q \) be a homogeneous prime in \( C \). Then \( \sigma _{Q}(f) \)
is a polynomial with coefficients in a field, hence \( b_{f_{Q}}(s) \) may
be computed. What would happen if we run Algorithm \ref{oakuAlg} trying to
compute \( b_{f_{Q}}(s) \) ``lifting from \( k(Q) \), the fraction field
of \( C/Q \), to \( C \)'' every single step of the algorithm? Notice that
\( \sigma _{Q}:C\rightarrow k(Q) \) has \( C/Q \) as its image. Since the
steps of the algorithm that do not involve Gröbner bases computation do not
involve division either, we have to worry only about the two steps that deal
with Gröbner bases. Suppose for these two steps we used \ref{alg_param} with
\( F' \) is a set generating \( Q \), in particular we obtained the exceptional
polynomials \( h_{1} \) and \( h_{2} \). Set \( h=h_{1}h_{2} \), then the
output, which is going to be \( b_{\sigma _{Q}(f)}(s) \), is also the Bernstein-Sato
polynomial of \( \sigma _{P}(f) \) for every \( P\supset Q \) such that \( h\notin P \).
Thus we have 

\begin{algorithm}
\label{1step} 
\end{algorithm}
\vspace{0.3cm}
{\centering \begin{tabular}{lll}
Input:&
\( f \):&
a polynomial in \( R_{n}(C) \),\\
&
\( F' \):&
generators of a homogeneous prime ideal,\\
Output:&
\( b(s) \):&
a polynomial in \( \Bbb {Q}[s] \), \\
&
\( H \):&
generators of a homogeneous ideal in \( C \) such that \\
&
&
\( b(s)=b_{\sigma _{P}(f)}(s) \) for every point \( P\in V'\backslash V'' \),
\\
&
&
where \( V'=V(F') \) and \( V''=V(H) \) (\( V''\subset V'\subset \Bbb {P}(n,d;k) \)).\\
\end{tabular}\par}
\vspace{0.3cm}

\begin{enumerate}
\item Compute the polynomial \( b(s) \) and the exceptional polynomial \( h \) as
described above. 
\item Return \( b(s) \) and \( \{h\}\cup F' \).
\end{enumerate}
\begin{rem}
\label{rem_tensor_k}If we consider \( C'=C\otimes k' \) and \( f\otimes 1\in R_{n}(C') \),
where \( k' \) is an extension of \( k \), then \( b(s) \) is the Bernstein-Sato
polynomial for any point in \( (V'\otimes _{k}k')\setminus (V''\otimes _{k}k') \).
\end{rem}
The next theorem gives an affirmative answer to Lyubeznik's question about the
constructibility of the set \( S(b(s),k) \) of Definition \ref{SBSdef}. 

\begin{thm}
The set \( S(b(s),k) \) is constructible for every \( b(s) \).
\end{thm}
\begin{proof}
The proof follows from the above algorithm. For the function \( \phi :\Bbb {P}(n,d;k)\rightarrow B(n,d) \),
\( \phi (P)=b_{\sigma _{P}(f)}(s) \) the following is true. For every projective
\( V'\subset \Bbb {P}(n,d;k) \) there is an open set \( U=V'\setminus V''\subset V' \)
such that \( \left. f\right| _{U} \) is a constant function. Therefore we may
apply Lemma \ref{lemma_constructible}.
\end{proof}
Algorithm \ref{1step} leads to the main algorithm and theorem of the paper. 

\begin{algorithm}
Input: \( n,d\in \Bbb {N} \). 

Output: The set of pairs \( L=\{(b(s),S(b(s)))|\, b(s)\in B(n,d)\} \), where
\( S(b(s))=S(b(s),\Bbb {Q})\subset \Bbb {P}(n,d;\Bbb {Q}) \).
\end{algorithm}
\begin{lyxlist}{00.00.0000}
\item [1.]Set \( L:=\emptyset  \), \( f:=\sum _{|\alpha |\leq d}a_{\alpha }x^{\alpha } \)
.
\item [2.]Define the recursive procedure \textbf{BSP}(\( Q \)), where \( Q\in \text {Spec}(\Bbb {Q}[\bar{a}]) \).
\end{lyxlist}
\vspace{0.3cm}
{\centering \begin{tabular}{l}
\textbf{BSP}(\( Q \)) := \{\\
\( \,  \)Apply Algorithm \ref{1step} to \( V(Q) \) and \( f \) \\
\hfill{}to get an ideal \( I \) in \( C \) and \( b(s)\in \Bbb {Q}[s] \);
\\
\( \,  \)IF there is a pair \( (b(s),S)\in L \) \\
\( \, \, \,  \)THEN replace it by \( (b(s),S\cup (V(Q)\setminus V(I))) \)
\\
\( \, \, \,  \)ELSE \( L:=L\cup \{(b(s),V(Q)\setminus V(I))\} \);\\
\( \,  \)IF \( V(I)\neq \emptyset  \) THEN \{ \\
\( \, \, \,  \)Find the minimal primes \( \{Q_{i}\} \) associated to \( I \);\\
\( \, \, \,  \)FOR each \( Q_{i} \) DO \textbf{BSP}(\( Q_{i} \)) ;\\
\( \,  \)\}\\
\}\\
\end{tabular}\par}
\vspace{0.3cm}

\begin{lyxlist}{00.00.0000}
\item [3.]Run BSP(\( 0 \)).
\end{lyxlist}
\begin{rem*}
This algorithm returns some presentations for constructible sets \( S(b(s),\Bbb {Q}) \),
the canonical presentations for which may be obtained by using the algorithm
discussed in Remark \ref{remark_constr_presentation}.
\end{rem*}
\begin{cor}
The set \( S(b(s),k) \) is defined over \( \Bbb {Q} \), i.e. there exist ideals
\( I_{i}\subset \Bbb {Q}[\bar{a}] \) and \( J_{i}\subset \Bbb {Q}[\bar{a}] \)
\( (i=1,...,m) \) such that for any field \( k \) 
\[
S(b(s),k)=\bigcup _{i}(V_{i}'\setminus V_{i}''),\]
 where \( V_{i}'=V(k[\bar{a}]I_{i}) \) is the zero set of the extension of
\( I_{i} \) and \( V_{i}''=V(k[\bar{a}]J_{i}) \) is the zero set of the extension
of \( J_{i} \). 
\end{cor}
\begin{proof}
Follows from the algorithm and Remark \ref{rem_tensor_k}.

\hfill{}
\end{proof}
The annihilators \( \text {Ann}(f^{s}) \) are computed using Algorithm \ref{AnnFs}
and the same technique as in the algorithm above. The output is a set of pairs
\( \{(I_{i},V_{i})\} \), where \( I_{i} \) are the ideals in \( A_{n}(k)[\bar{a}][s] \)
and \( V_{i} \) are locally closed sets, such that for any polynomial \( f \)
with coefficients in \( k \) that corresponds to a point \( P\in V_{i} \)
the ideal \( \text {Ann}(f^{s}) \) equals \( \sigma _{P}(I_{i}) \), the ideal
\( I_{i} \) specialized to \( P \).

After doing the above steps, the real life algorithm that produces Example \ref{exaAnnFs}
compresses its output in the following way. If \( (I_{i},V_{i}) \) and \( (I_{j},V_{j}) \)
are two different pairs such that \( \sigma _{P}(I_{i})=\sigma _{P}(I_{j}) \)
for all \( P\in V_{j} \) then these two are replaced by the pair \( (I_{i},V_{i}\cup V_{j}) \). 

\begin{rem}
\label{rem_quasidefinitive}The stratification of the parameter space produced
by such computation is not unique. So we have to use prefix ``quasi`` in ``quasi-definitive
computation'', because the annihilators, as opposed to Bernstein-Sato polynomials,
depend on the parameters making it possible to slice the space of parameters
in many ways. 
\end{rem}

\end{document}